\documentclass{article}
\usepackage[top=30truemm,bottom=30truemm,left=25truemm,right=25truemm]{geometry}
\usepackage{amsmath,amssymb,amsfonts,mathtools,accents}
\usepackage[dvipdfmx]{graphicx}
\usepackage{psfrag,epsf}
\usepackage{enumerate}
\usepackage{tabu}
\newtheorem{theorem}{Theorem}

\newtheorem{cor}{Corollary}
\newtheorem{assumption}{Assumption}
\newtheorem{remark}{Remark}
\usepackage{url} 

\title{A new method of joint nonparametric estimation of probability density and its support}
\author{
Taku MORIYAMA\\
Graduate School of Mathematics, Kyushu University}
\date{}

\begin{document}
\maketitle

\begin{abstract}
We propose a new method for simultaneous nonparametric estimation of a probability density and its support. As is well known, a nonparametric kernel density estimator has a `boundary bias problem' when the support of the population density does not cover the whole real line. If we know the support exactly, we may reduce the bias by using a boundary bias reduction method. When the support is unknown, there is possibly a boundary problem of which we should take care. We insist on the necessity of estimating the support and propose a new method of nonparametric density estimation that is free from the boundary bias in such case.

The proposed method detects the boundary and gives the modified density estimator simultaneously. Although it is natural to estimate the support by using the sample maximum (and minimum) and modify the density estimator as in Hall \& Park (2002), the new method is numerically superior in the sense of an integrated squared error in the boundary region. Moreover, we discuss an extension to a simple multivariate case and propose a new method for estimating the joint probability density. Using the idea of nonparametric copula estimation, this method combines the marginal densities estimated by the proposed single variable method. The obtained joint density estimator is also boundary bias free.
\end{abstract}

\section{Introduction}
\label{sec:intro}
Random phenomena are described uniquely by their probability distribution, and estimations of their distributions give us much information. Nonparametric density estimation not only provides a graphical overview of the shape of a distribution, but also enables us to infer a variety of interesting things. These include estimation of functionals of density and statistical testing, for example, goodness of fit with a parametric model, equality of two distributions, symmetry, multi-modality, and so on. Rosenblatt (1956) has proposed a smooth nonparametric estimator $\widehat{f}$ of the density function $f$, and it has been extensively investigated (Tsybakov (2009)). The kernel cumulative distribution estimator is given by integration of $\widehat{f}$; we will begin by introducing that.

Let $X_1,X_2,\cdots,X_n$ be independently and identically distributed ($i.i.d.$) random variables with a distribution function $F$, and let $f$ be the density function. The kernel density estimator and distribution estimator are given by
$$\widehat{f}(x):=\frac{1}{nh} \sum_{i=1}^n K\left(\frac{x -X_i}{h} \right)$$
and
$$\widehat{F}(x):=\frac{1}{n} \sum_{i=1}^n W\left(\frac{x -X_i}{h} \right),$$
where $K$ is a symmetric kernel function, $W$ is the integral of $K$ and $h$ is a bandwidth which satisfies $h \rightarrow 0$ and $n h \rightarrow \infty ~(n \rightarrow \infty)$. We call $\widehat{f}$ the `naive density estimator' and $\widehat{F}$ similarly.

Since the naive kernel estimators are sums of $i.i.d.$ random variables, it is easy to obtain their moments. By making a change of variables and performing a Taylor expansion, if the support of $f$ is all of $\mathbb{R}$, we find that
$$E[\widehat{f}(x)] = \frac{1}{h} \int_{-\infty}^{\infty} K\left(\frac{x -y}{h} \right) f(y) = \int_{-\infty}^{\infty} K\left(z \right) f(x -hz) = f(x) + O(h^2)$$
and similarly
$$V[\widehat{f}(x)] = \frac{f(x)}{nh} + O\left(\frac{1}{n}\right),$$
$$E[\widehat{F}(x)] = F(x) + O(h^2), ~~~~~ V[\widehat{F}(x)] = \frac{F(x)}{n} + O\left(\frac{h}{n} \right)$$
under certain regularity conditions. However, when the support is not all of $\mathbb{R}$, the moments changes, and $\widehat{f}$ loses consistency near the boundary of the support. This situation is known as the `boundary bias' problem. When the support is $(-\infty,u_0]$ (i.e. $f(u_0)>0$), the bias of $\widehat{f}$ at the boundary is
\begin{eqnarray*}
Bias[\widehat{f}(u_0)] &=& \frac{1}{h} \int_{-\infty}^{u_0} K\left(\frac{u_0 -y}{h} \right) f(y) dy -f(u_0)\\
&=& \int_{0}^{\infty} K\left(z \right) f(u_0 -hz) dz -f(u_0) = -\frac{f(u_0)}{2} + O(h).
\end{eqnarray*}

We call the area $(u_0-h,u_0]$, where $\widehat{f}$ is biased, the boundary region of the density $f$. In addition, the order of the asymptotic bias of $\widehat{F}$ becomes larger: for $x \in \mathbb{R}$,
$$Bias[\widehat{F}(x)] =O(h) ~~~\text{uniformly}.$$

\begin{remark}
In fact, the boundary bias depends on the left limit $f(u_0-)$ rather than $f(u_0)$. Therefore, it is not a problem of whether $f$ is left or right continuous at $u_0$. Hereafter, we will assume that $f(u_0-)=f(u_0)$ and that $f$ is right continuous at the left endpoint $l_0$ $($i.e. $f(l_0+)=f(l_0)$$)$ if $l_0$ is bounded.
\end{remark}

If we know the support of $f$ exactly, we can reduce the bias by using a boundary bias reduction method. The basic reduction methods of $\widehat{f}$ include renormalization (Jones (1993)), reflection (\'{C}wik \& Mielniczuk (1993)), asymmetric kernel (Chen (1999), Chen (2000), {\it etc}.) and generalized jackknife (Jones (1993), Terrell \& Scott(1980)), the `direct' method for reducing to an arbitrary order (Bearse \& Rilstone (2009)), and so on. It seems to be easy to obtain boundary bias free estimators of the cumulative distribution $F$ by integration of the modified kernel density estimators. However, the modified density estimators are not usually $i.i.d.$ sums, and it is often difficult to represent the integral explicitly. That is why it is hard obtain the modified cumulative distribution estimators. Recently, some papers were published that focus on distribution function estimation. The boundary kernel method (Tenreiro(2013)) and generalized reflection method (Kolacek \& Karunamuni (2011)) can give boundary bias free estimators of $F$. We should note that all of them control the boundary effect which comes from the `known' support.

When the support is unknown, however, we do not take care with the boundary problem well. This may be because there are almost no papers tackling the unknown boundary effect (except Hall \& Park(2002)). In most `actual' cases, the realized values are not large, and so the support must be smaller than $\mathbb{R}$. Therefore, we insist on the necessity of both estimating the support and eliminating the unknown boundary bias appropriately.

To overcome the unknown boundary effect, one can estimate the support and then modify the density estimator which regards the estimated support as true. In fact, Hall \& Park(2002) proposed to replace the unknown upper bound of the support with the sample maximum $X_{(n)}$. We call the modified density estimator an `$X_{(n)}$-based estimator'. However, it seems obvious that the boundary estimator $X_{(n)}$ is not always numerically best for a density estimation that utilizes a very different boundary bias reduction method. Here, we propose a new method for estimating the probability density and its support simultaneously. The boundary estimator depends on the boundary bias reduction method, which is what we apply, and the proposed density estimator is numerically more accurate in the boundary region. In fact, the new method minimizes a loss function asymptotically. 

Section 2 describes some of the basic boundary bias reduction methods for when the support is known and simple, such as $[0,\infty)$ and $[0,1]$. In Section 3, we describe the new method for estimating the population density $f$ which is free from the `unknown' boundary bias and investigate its asymptotic properties. In addition, we confirm that some boundary bias reduction methods satisfy a condition which ensures the proposed estimator works well. In Section 4, we compare the proposed estimator with the naive kernel density estimator $\widehat{f}$ and the $X_{(n)}$-based estimator in the sense of the integrated squared error in the boundary region. When $u_0$ is bounded and $f(u_0)>0$ holds, $X_{(n)}$ has $n$ consistency, while the new boundary estimator has $\sqrt{n}$ consistency. However, the convergence rate does not affect the integrated error asymptotically; in fact, we demonstrate that the new density estimator performs better numerically. Moreover, in section 5, we discuss an extension to the simple multivariate case in which the support is given by an hyper-rectangle. Using the idea of nonparametric copula estimation, we propose a new method for estimating joint density estimator. This method combines the marginals estimated by the single variable method (one-dimensional cases), and the obtained joint density estimator is boundary bias free. We study the proposed method by simulating a number of bivariate distribution estimations. The proofs are given in the appendices.

\section{Boundary bias reduction methods}
\label{Boundary bias reduction methods}
Let us assume that the support of the density function $f$ is known and given by $[l_0,u_0]$. \'{C}wik \& Mielniczuk (1993) discussed a `reflection' method that reduces the boundary bias of the naive kernel density estimator. The estimator is given by
\begin{align*}
&\widehat{f}^{[R]}(x) := \widehat{f}(x) + \widehat{f}(2u_0 -x) + \widehat{f}(2l_0 -x) \\
=&\frac{1}{nh} \sum_{i=1}^n \biggl \{ K\left(\frac{x -X_i}{h} \right) +K \left(\frac{x +X_i -2u_0}{h} \right) +K \left(\frac{x +X_i -2l_0}{h} \right) \biggl \}
\end{align*}
for $x \in [l_0,u_0]$. The asymptotic bias is as follows: for $x \in \mathbb{R}$,
$$Bias[\widehat{f}^{[R]}(x)] =O(h) =:b_f^{[R]}(x) h \hspace{10pt}\text{(uniformly)},$$
where $b_f^{[R]}(x)=O(1)$, and the asymptotic order of the variance is the same as that of the naive estimator $\widehat{f}$. Thus, the estimator recovers its consistency in the boundary region. In addition, if $f'(u_0-)=0$, the bias becomes of order $h^2$ uniformly, where $f'$ is the derivative of $f$. From the integral of $\widehat{f}^{[R]}(x)$, we can derive the distribution estimator as follows:
\begin{align*}
&\widehat{F}^{[R]}(x) := \int_{l_0}^{x} \widehat{f}^{[R]}(y) dy\\
=& \frac{1}{n} \sum_{i=1}^n \bigg \{ W\left(\frac{x -X_i}{h} \right) +W \left(\frac{x +X_i -2u_0}{h} \right) +W \left(\frac{x +X_i -2l_0}{h} \right) \\
& ~~~~~~ -W \left(\frac{X_i +l_0 -2u_0}{h} \right) \bigg \} -1 \\
=& \widehat{F}(x) -\widehat{F}(2u_0 -x) -\widehat{F}(2l_0 -x) + \widehat{F}(2u_0 -l_0) ~~~ (x \in [l_0,u_0]).
\end{align*}
From \'{C}wik \& Mielniczuk (1993, we can see
$$Bias[\widehat{F}^{[R]}(x)] =O(h^2) =:b_F^{[R]}(x) h^2 \hspace{10pt}\text{(uniformly for $x \in \mathbb{R}$)},$$
where $b_F^{[R]}(x)=O(1)$, and the order of the variance is the same as that of the naive estimator.

Tenreiro (2013) proposed the `boundary kernel' method for reducing the bias of the kernel cumulative distribution estimator. The method changes the kernel function in the boundary region. It has a simple form, given as follows:

\[
  \widehat{F}^{[BK]}(x) = \begin{dcases}
    \widehat{F}^{[BK,l_0]}(x) := \frac{1}{n} \sum_{i=1}^n W \left(\dfrac{x -X_i}{x -l_0} \right) &(l_0 \le x <l_0+h) \\[2ex]
    \widehat{F}^{[BK,I]}(x) := \widehat{F}(x) &(l_0+h \le x < u-h) \\[2ex]
    \widehat{F}^{[BK,u_0]}(x) := 1 - \frac{1}{n} \sum_{i=1}^n W \left(\dfrac{X_i -x}{u_0 -x} \right) &(u_0-h \le x < u_0) \\[2ex]
    0, ~~~ 1 &(x < l_0), ~~~(u_0 \le x).
  \end{dcases}
\]
The author of that paper shows that for $x \in \mathbb{R}$,
$$Bias[\widehat{F}^{[BK]}(x)] =O(h^2) =: b_F^{[BK]}(x) h^2 \hspace{10pt}\text{(uniformly)}$$
when the support of $K$ is bounded. We can also derive the following density estimator $\widehat{f}^{[BK]}$ from the derivative of $\widehat{F}^{[BK]}$ and calculate the following bias,
$$Bias[\widehat{f}^{[BK]}(x)] =O(h) =: b_f^{[BK]}(x) h \hspace{10pt}\text{(uniformly)}.$$
If $f'(u_0-)$ is zero, the bias becomes of order $h^2$.

\section{New joint estimator of probability density and its support in one dimension}
\label{New joint estimator of probability density and its support in 1-dimension} 
Let us assume that both boundaries $\mbox{\boldmath $u$}_0 =(l_0,u_0)^T$ are unknown but that both are bounded. We do not assume that the interval is open, half-open, or closed. Now we put $\mbox{\boldmath $u$} =(l,u)^T$, and $\widehat{f}_{\mbox{\boldmath $u$}}^{\dagger}$ denotes the kernel type and boundary bias free estimator of the density whose support is $[l,u]$. We propose the following new method for estimating the probability density and its support simultaneously. Let us define the estimator $\mbox{\boldmath $u$}=\widehat{\mbox{\boldmath $u$}}$ of $\mbox{\boldmath $u$}_0$ as the solution of the following function,
\begin{eqnarray}
\widehat{\mbox{\boldmath $F$}}_{\mbox{\boldmath $u$}}^{\dagger} ( \mbox{\boldmath $X$}_{(1,n)} ) - \mbox{\boldmath $c$}_n = \mbox{\boldmath $0$}
\end{eqnarray}
and the density estimator as $\widehat{f}_{\widehat{\mbox{\boldmath $u$}}}^{\dagger}$ (and the distribution estimator as $\widehat{F}_{\widehat{\mbox{\boldmath $u$}}}^{\dagger}$), where
\begin{eqnarray}
\widehat{\mbox{\boldmath $F$}}_{\mbox{\boldmath $u$}}^{\dagger}((x_1, x_2)^T) := \left( \widehat{F}_{\mbox{\boldmath $u$}}^{\dagger}(x_1) ,   \widehat{F}_{\mbox{\boldmath $u$}}^{\dagger}(x_2) \right)^T = \left( \int_{-\infty}^{x_1}  \widehat{f}_{\mbox{\boldmath $u$}}^{\dagger}(z) dz ~,~ \int_{-\infty}^{x_2}  \widehat{f}_{\mbox{\boldmath $u$}}^{\dagger}(z) dz \right)^T
\end{eqnarray}
and
\begin{eqnarray}
\mbox{\boldmath $X$}_{(1,n)} = \left(X_{(1)}, X_{(n)} \right)^T
,~~~
\mbox{\boldmath $c$}_{n} = 
\left(\frac{1}{n+1}~,~ \frac{n}{n+1} \right)^T
,~~~ 
\mbox{\boldmath $0$}=(0,0)^T.
\end{eqnarray}
In this section, all functions in bold-type denote two-dimensional vectors of `same' scalar-valued functions for  (e.g., $\mbox{\boldmath $F$}(\mbox{\boldmath $x$})=(F(x_1),F(x_2))^T$). Note that $\mbox{\boldmath $u$}_0$ stands for the true value and $\mbox{\boldmath $u$}$ is variable. Intuitively, equation $(1)$ comes from the properties of the maximum and minimum estimation of the uniform distribution on the interval $[0,1]$. This is because $\mbox{\boldmath $F$}((X_1, X_2)^T)  \stackrel{d}=\mbox{\boldmath $Z$}$ (in distribution) and $\mbox{\boldmath $F$}(\mbox{\boldmath $X$}_{(1,n)})\approx \mbox{\boldmath $c$}_n$ hold, where $\mbox{\boldmath $Z$}$ is the two-dimensional uniform random variable on $[0,1] \times [0,1]$. Moreover, $(1+(1/n))Z_{(n)}$ (the sample maximum of the uniform distribution) is known as the minimum-variance unbiased estimator of the boundary value ($u_0=1$) of the uniform distribution (the sample minimum is similar).

\begin{remark}
The solution $\mbox{\boldmath $u$}=\mbox{\boldmath $X$}_{(1,n)}$ (that is, $\mbox{\boldmath $X$}_{(1,n)}$-based estimator) can be viewed as the solution of the following equation,
$$\widehat{\mbox{\boldmath $F$}}_{{\mbox{\boldmath $u$}}}^{\dagger}(\mbox{\boldmath $X$}_{(1,n)}) - \mbox{\boldmath $c$} =\mbox{\boldmath $0$},$$
where $\mbox{\boldmath $c$}=(0,1)^T$.
\end{remark}

In fact, the new estimator asymptotically coincides with the minimizer of a local expected predict error. Let us define the following error as the loss function
\begin{align*}
&E_{X_{n+1}}\left[\Bigm(\widehat{F}_{\mbox{\boldmath $u$}}^{\dagger}(X_{n+1}) - F(X_{n+1}) \Bigm)^2 I(X_{n+1} \in \mbox{\boldmath $q$}_n) \Bigm| X_1, \cdots, X_{n}\right] \\
=& \int_{\mbox{\boldmath $q$}_n} \Bigm(\widehat{F}_{\mbox{\boldmath $u$}}^{\dagger}(x) - F(x) \Bigm)^2 f(x) dx,
\end{align*}
where $\widehat{F}_{\mbox{\boldmath $u$}}^{\dagger}$ is based on the sample $(X_1, \cdots, X_n)$ and $I$ is the indicator function $I(A)=1$ (if $A$ occurs), $=0$ (if $A$ fails). Let $\mbox{\boldmath $q$}_n$ be some two areas near both boundaries which contain both $X_{(1)}$ and $X_{(n)}$ as the realized value. Then, replacing $f(x) dx$ by the estimator $dF_n(x)$ and using the fact that $\mbox{\boldmath $F$}(\mbox{\boldmath $X$}_{(1,n)})\approx \mbox{\boldmath $c$}_n$, we can see that the minimizer of the error asymptotically coincides with $\widehat{\mbox{\boldmath $u$}}$ as follows:
\begin{align*}
0 &= \int_{\mbox{\boldmath $q$}_n} \Bigm(\widehat{F}_{\mbox{\boldmath $u$}}^{\dagger}(x) - F(x) \Bigm)^2 dF_{n}(x) \\
\approx& \Bigm\|\widehat{\mbox{\boldmath $F$}}_{\mbox{\boldmath $u$}}^{\dagger}(\mbox{\boldmath $X$}_{(1,n)}) - \mbox{\boldmath $c$}_n \Bigm\|^2.
\end{align*}

The solution $\mbox{\boldmath $u$}= \widehat{\mbox{\boldmath $u$}}$ of $(1)$ is not given as an explicit formula, and the properties depend on the bias reduction method. Next, we state the general properties of the proposed estimator and study some applications.

\subsection{Asymptotic properties}
To construct the asymptotic properties of the new estimators, we utilize the asymptotic theory of $M$ estimation. To view the solution $\mbox{\boldmath $u$}= \widehat{\mbox{\boldmath $u$}}$ as an $M$ estimator asymptotically, we make the following assumptions.

\begin{assumption}
For all large enough integers $n$, there is $\mbox{\boldmath $\Theta$}$ whose interior includes $\mbox{\boldmath $u$}_0$. 
\end{assumption}

\begin{assumption}
There is a function $\mbox{\boldmath $\Psi$}_{\mbox{\boldmath $u$},n}$ which satisfies the following  for all large enough integers $n$: \\
$(i)$ $\mbox{\boldmath $\Psi$}_{\mbox{\boldmath $u$},n}$ is given by following $i.i.d.$ sum
$$\mbox{\boldmath $\Psi$}_{\mbox{\boldmath $u$},n}(\mbox{\boldmath $x$})  = \frac{1}{n} \sum_{i=1}^n \mbox{\boldmath $\psi$}_{i,\mbox{\boldmath $u$},n}(\mbox{\boldmath $x$})$$
where
$$\mbox{\boldmath $\Psi$}_{\mbox{\boldmath $u$},n}(\mbox{\boldmath $x$}) := (\Psi_{\mbox{\boldmath $u$},n}(x_1), \Psi_{\mbox{\boldmath $u$},n}(x_2))^T$$
and 
$$\mbox{\boldmath $\psi$}_{i,\mbox{\boldmath $u$},n}(\mbox{\boldmath $x$}) := (\psi_{i,\mbox{\boldmath $u$},n}(x_1), \psi_{i,\mbox{\boldmath $u$},n}(x_2))^T.$$
$(ii)$ The following stochastic expansion holds:
\begin{eqnarray}
\widehat{\mbox{\boldmath $F$}}_{\mbox{\boldmath $u$}}^{\dagger}(\mbox{\boldmath $x$}) = \mbox{\boldmath $\Psi$}_{\mbox{\boldmath $u$},n}(\mbox{\boldmath $x$}) + \mbox{\boldmath $R$}_{\mbox{\boldmath $u$},n}(\mbox{\boldmath $x$}),
\end{eqnarray}
~ where the residual $\mbox{\boldmath $R$}_{\mbox{\boldmath $u$},n}$ satisfies
$$\mbox{\boldmath $R$}_{\mbox{\boldmath $u$},n} (\mbox{\boldmath $x$})=o_P(n^{-1/2} \mbox{\boldmath $1$}),~~~ \frac{\partial}{\partial \mbox{\boldmath $u$}} \mbox{\boldmath $R$}_{\mbox{\boldmath $u$},n}(\mbox{\boldmath $x$}) =o_P(n^{-1/2} \mbox{\boldmath $1$})$$
~ uniformly for $\mbox{\boldmath $u$} \in \mbox{\boldmath $\Theta$}$ and $\mbox{\boldmath $x$} \in \mathbb{R}^2$, and $\mbox{\boldmath $1$}=(1,1)^T$. \\
$(iii)$ The next convergence holds
\begin{eqnarray}
\sup_{\mbox{\boldmath $u$} \in \mbox{\boldmath $\Theta$}} \Biggm\|\mbox{\boldmath $\Psi$}_{\mbox{\boldmath $u$},n}(\mbox{\boldmath $u$}_0) -\mbox{\boldmath $F$}_{\mbox{\boldmath $u$}}(\mbox{\boldmath $u$}_0) \Biggm\| \xrightarrow{p} 0,
\end{eqnarray}
~ where $\mbox{\boldmath $F$}_{\mbox{\boldmath $u$}}$ is a function such that $\mbox{\boldmath $F$}_{\mbox{\boldmath $u$}_0}$ equals the  underlying \\
~ distribution function $F$.
\end{assumption}

\begin{assumption}
The following holds for all large enough integers $n$: \vspace{1ex}\\
$(i)$ $\centering{ \displaystyle \Psi_{\mbox{\boldmath $u$}_0,n}(x) = F(x) + b_F^{\dagger}(x,n,h) + O_P(n^{-1/2}) }$ \vspace{1ex} \\
~ holds uniformly for $x \in \mathbb{R}$, where $b_F^{\dagger}$ satisfies $b_F^{\dagger}(\mbox{\boldmath $u$}_0,n,h) = o(n^{-1/2} \mbox{\boldmath $1$})$.\vspace{1ex} \\
$(ii)$ $\displaystyle \frac{\partial}{\partial x}\Psi_{\mbox{\boldmath $u$}_0,n}(x) = f(x) + b_f^{\dagger}(x,n,h) + O_P((nh)^{-1/2})$ \vspace{1ex}\\
~ holds uniformly for $x \in \mathbb{R}$, where
$$\frac{\partial}{\partial x}\Psi_{\mbox{\boldmath $u$}_0,n}(x_0) = \frac{\partial}{\partial x}\Psi_{\mbox{\boldmath $u$}_0,n}(x) \biggm|_{x=x_0}.$$
$(iii)$ For every $l~(l=1,2)$,
\begin{eqnarray*}
\left(\bigotimes_{k=1}^l \frac{\partial}{\partial \mbox{\boldmath $u$}}\right) \Psi_n(x, \mbox{\boldmath $u$}) = O\left(\bigotimes_{k=1}^l \mbox{\boldmath $1$}\right) + o_P\left(\bigotimes_{k=1}^l \mbox{\boldmath $1$}\right)
\end{eqnarray*}
~ holds uniformly for $\mbox{\boldmath $u$} \in \mbox{\boldmath $\Theta$}$ and $x \in \mathbb{R}$ (that is, every component converges some uniformly bounded constants), where
\begin{eqnarray*}
\frac{\partial}{\partial \mbox{\boldmath $u$}} \Psi_{\mbox{\boldmath $u$}_0,n}(x) = \frac{\partial}{\partial \mbox{\boldmath $u$}} \Psi_{\mbox{\boldmath $u$},n}(x)  \biggm|_{\mbox{\boldmath $u$}=\mbox{\boldmath $u$}_0}.
\end{eqnarray*}
\end{assumption}

\begin{assumption}
For all large enough integers $n$, the following holds: \vspace{1ex}\\
$(i)$ ~ $E[\mbox{\boldmath $\psi$}_{i,\mbox{\boldmath $u$},n}(\mbox{\boldmath $u$}_0)] -\mbox{\boldmath $c$}_n = \mbox{\boldmath $0$}$ ~\text{has the unique solution}~ $\mbox{\boldmath $u$}=\mbox{\boldmath $u$}_n^*$ \text{which satisfies}~ \vspace{1ex}\\
 ~~~ $\mbox{\boldmath $u$}_n^* \to \mbox{\boldmath $u$}_0$. \vspace{1ex}\\
$(ii)$ ~ \text{for any $\eta>0$, there is $\kappa_{\eta,\delta_n,n} = O(1)$ which ensures} \vspace{1ex}\\
~~~ ${\displaystyle P\left[\inf_{\mbox{\boldmath $u$}: \|\mbox{\boldmath $u$} -\mbox{\boldmath $u$}_0\|> \eta} \| \mbox{\boldmath $\Psi$}_{\mbox{\boldmath $u$},n}(\mbox{\boldmath $u$}_0) -\mbox{\boldmath $c$}_n\|  > \kappa_{\eta,\delta_n,n} \right] > 1-\delta_n }$, \vspace{1ex}\\
 ~~~ \text{where $\delta_n=o(1) >0$}. \vspace{1ex}\\
$(iii)$ ~ ${\displaystyle \left(E\left[\frac{\partial}{\partial \mbox{\boldmath $u$}} \mbox{\boldmath $\psi$}_{i,\mbox{\boldmath $u$}_n^*,n}(X_i,\mbox{\boldmath $u$}_0) \right]\right)^{-1} }$ ~~\text{exists, i.e., the matrix is nonsingular.}\vspace{1ex}\\
$(iv)$ ~ $E[ \|\mbox{\boldmath $\psi$}_{i,\mbox{\boldmath $u$},n}(\mbox{\boldmath $u$}_0)\|^2]$ ~\text{is bounded in some neighborhood of}~ $\mbox{\boldmath $u$}=\mbox{\boldmath $u$}_n^*$\vspace{1ex}\\
~\text{and is continuous at}~ $\mbox{\boldmath $u$}_n^*$.
\end{assumption}

Assumptions 2 and 3 admit the following asymptotic expansion of $\widehat{\mbox{\boldmath $F$}}_{\mbox{\boldmath $u$}}^{\dagger}(\mbox{\boldmath $X$}_{(1,n)})$ at the $i.i.d.$ sum $\mbox{\boldmath $\Psi$}_{\mbox{\boldmath $u$},n}$, 
\begin{eqnarray}
\widehat{\mbox{\boldmath $F$}}_{\mbox{\boldmath $u$}}^{\dagger}(\mbox{\boldmath $X$}_{(1,n)}) = \mbox{\boldmath $\Psi$}_{\mbox{\boldmath $u$},n}(\mbox{\boldmath $u$}_0) + o_P(n^{-1/2}\mbox{\boldmath $1$}),
\end{eqnarray}
and the boundary bias reduction of the distribution estimators $\widehat{F}_{\mbox{\boldmath $u$}}^{\dagger}$ and $\widehat{f}_{\mbox{\boldmath $u$}}^{\dagger}$. Assumption 4 keeps the $\sqrt{n}$ consistency of $\widehat{\mbox{\boldmath $u$}}$, on the basis of the asymptotic theory of $M$ estimation.

Under these assumptions, we can show the consistency of $\widehat{\mbox{\boldmath $u$}}$ and the bias reduction of $\widehat{F}_{\widehat{\mbox{\boldmath $u$}}}^{\dagger}$ and $\widehat{f}_{\widehat{\mbox{\boldmath $u$}}}^{\dagger}$.

\begin{theorem} 
Given Assumptions 1 - 4 and that
\begin{eqnarray}
\sqrt{n}(\mbox{\boldmath $X$}_{(1,n)} -\mbox{\boldmath $u$}_0) \xrightarrow{p} \mbox{\boldmath $0$}
\end{eqnarray}
we have
\begin{align}
\widehat{\mbox{\boldmath $u$}} =& \mbox{\boldmath $u$}_0 + O_P(n^{-1/2} \mbox{\boldmath $1$}), \\
\widehat{F}_{\widehat{\mbox{\boldmath $u$}}}^{\dagger} (x) =& F(x) + b_F^{\dagger}(x,n,h) + O_P(n^{-1/2}) \\ 
\intertext{and}
\widehat{f}_{\widehat{\mbox{\boldmath $u$}}}^{\dagger} (x) =& f(x) + b_f^{\dagger}(x,n,h) + O_P((nh)^{-1/2}).
\end{align}
\end{theorem}
{\it Proof.} See the Appendices.

\begin{remark}
When $f(l_0)>0$ and $f(u_0)>0$, $\mbox{\boldmath $X$}_{(1,n)}$ is an $n$ consistent estimator of $\mbox{\boldmath $u$}_0$ and clearly satisfies the assumptions of Theorem 1. 
\end{remark}

\begin{remark}
If we know either $l_0$ or $u_0$, we can see that Theorem 1 still holds by replacing the all vector values with the scalar ones, for example, $\mbox{\boldmath $u$}_0$ with a scalar $u_0$ or $l_0$, $\mbox{\boldmath $u$}$ with a scalar $u$ or $l$, and $\mbox{\boldmath $X$}_{(1,n)}$ with $X_{(n)}$ or $X_{(1)}$, and so on. Table 1 shows the main procedure of the proposed method if we know the left support $l_0$.
\end{remark}

\begin{table}
\caption{Flow of the proposed method (the left boundary $l_0$ is known)}
\begin{tabular}{l}
\hline \vspace{3pt}
0. Assume that the support of the density is from $l_0$ to $u_0$, \\ \vspace{3pt}
~~~~~ where $l_0$ is known, and $u_0$ is unknown but bounded. \\ \vspace{3pt}
~~~ Check the assumptions of Theorem 1. \\ \vspace{3pt}
1. Select a boundary bias reduction method, \\ \vspace{3pt}
~~~~~ and choose the kernel and bandwidth, following the method.\\ \vspace{3pt}
~~~ $\widehat{f}_{u}^{\dagger}$ denotes the boundary bias free estimator \\
~~~~~ of the density function whose support is $l_0$ to $u$ \\
~~~~~ and $\widehat{F}_{u}^{\dagger}$ denotes the cumulative distribution estimator.\\
\vspace{5pt}
2. Solve the equation for $u$ : \\ \vspace{5pt}
\hspace{100pt} $\widehat{F}_{u}^{\dagger}(X_{(n)}) - c_n =0$. \\ \vspace{3pt}
3. Set the solution as $\widehat{u}$, \\ \vspace{3pt}
~~~~~ and output the boundary-adjusted estimators $\widehat{F}_{\widehat{u}}^{\dagger}$ and $\widehat{f}_{\widehat{u}}^{\dagger}$.\\
\hline
\end{tabular}
\end{table}

Now, we apply the boundary kernel method and reflection method described above. Since $P[\|\mbox{\boldmath $X$}_{(1,n)} -\mbox{\boldmath $u$}_0\|> h] = o((\sqrt{n}h)^{-1})$ holds under the assumptions of Theorem 1, we define the boundary estimator $\mbox{\boldmath $u$}=\widehat{\mbox{\boldmath $u$}}^{[BK]}$ as the solution of the following equation,
\begin{eqnarray}
\widehat{\mbox{\boldmath $F$}}_{\mbox{\boldmath $u$}}^{[BK,\mbox{\boldmath $u$}]}(\mbox{\boldmath $X$}_{(1,n)}) - \mbox{\boldmath $c$}_{n} = \mbox{\boldmath $0$},
\end{eqnarray}
where 
\begin{eqnarray*}
\widehat{\mbox{\boldmath $F$}}_{\mbox{\boldmath $u$}}^{[BK,\mbox{\boldmath $u$}]}(\mbox{\boldmath $X$}_{(1,n)}) = \left(\widehat{F}_{l}^{[BK,l]}(X_{(1)}),~ \widehat{F}_{u}^{[BK,u]}(X_{(n)}) \right)^T.
\end{eqnarray*}
The vector-valued equation $(11)$ is divided into two scalar-valued equations, such that one includes only $l$, and the other includes only $u$. We can easily see that $\widehat{F}_{l}^{[BK,l]}(X_{(1)})$ is an increasing function for $l \in (-\infty,X_{(1)}]$ and $\widehat{F}_{u}^{[BK,u]}(X_{(n)})$ is a decreasing function for $u \in [X_{(n)},\infty)$. For any fixed $n (\ge 2)$, we have
\begin{align*}
\lim_{u \downarrow X_{(n)}} \widehat{F}_{u}^{[BK,u]}(X_{(n)}) &= 1- \frac{1}{2n} ~~~ \left(\ge \frac{n}{n+1} \right)
\end{align*}
and
$$\lim_{u \uparrow \infty} \widehat{F}_{u}^{[BK,u]}(X_{(n)}) = \frac{1}{2} ~~~ \left(\le \frac{n}{n+1} \right).$$
Therefore, we can see that the solution $u=\widehat{u}$ is unique and $l=\widehat{l}$ is also unique. Thus, it is easy to confirm that the minimizer $\widehat{\mbox{\boldmath $u$}}^{[BK]}$ is unique and satisfies $(11)$. Now, we can deduce the following result.

\begin{cor} 
Let us assume that $supp(K)=[-1,1]$ and $f^{(1)}$ exists and is continuous. When $(7)$ holds and $h=o(n^{-1/4})$, we have
\begin{align*}
\widehat{\mbox{\boldmath $u$}}^{[BK]} =& \mbox{\boldmath $u$}_0 + O_P(n^{-1/2} \mbox{\boldmath $1$}),\\
\widehat{F}_{\widehat{\mbox{\boldmath $u$}}}^{[BK]} (x) =& F(x) + b_F^{[BK]}(x) h^2 + O_P(n^{-1/2})\\
\intertext{and}
\widehat{f}_{\widehat{\mbox{\boldmath $u$}}}^{[BK]} (x) =& f(x) + b_f^{[BK]}(x) h + O_P((nh)^{-1/2}).
\end{align*}
\end{cor}

By applying the reflection method, we define the boundary estimator $\mbox{\boldmath $u$}=\widehat{\mbox{\boldmath $u$}}^{[R]}$ and the distribution estimator $\widehat{F}_{\widehat{\mbox{\boldmath $u$}}}^{[R]}$ in the same way:
\begin{eqnarray}
\widehat{\mbox{\boldmath $F$}}_{\mbox{\boldmath $u$}}^{[R]}(\mbox{\boldmath $X$}_{(1,n)}) -\mbox{\boldmath $c$}_{n} = \mbox{\boldmath $0$}.
\end{eqnarray}
When $n$ is large enough, equation $(12)$ is divided into two mutually independent parts in the sense of $l$ and $u$ as same as the boundary kernel method. It is easy to see that $\widehat{F}_{l}^{[R]}(X_{(1)})$ is an increasing function for $l \in (-\infty,X_{(1)}]$ and that $\widehat{F}_{u}^{[R]}(X_{(n)})$ is a decreasing function for $u \in [X_{(n)},\infty)$. For any fixed $n$, we have
\begin{align*}
\lim_{u \downarrow X_{(n)}} \widehat{F}_{u}^{[R]}(X_{(n)}) &= \frac{1}{n} \sum_{i=1}^n \left\{W\left( \frac{X_{(n)} +X_i -2l_0}{h}\right) - W\left( \frac{ l_0 -X_{(n)} }{h}\right) \right\} (=: \widehat{F}_{X_{(n)}}^{[R]}(X_{(n)}))
\end{align*}
and
\begin{align*}
\lim_{u \uparrow \infty} \widehat{F}_{u}^{[R]}(X_{(n)}) &= \frac{1}{n} \sum_{i=1}^n \left\{ W\left( \frac{X_{(n)} +X_i -2l_0}{h}\right) - W\left( \frac{X_i -X_{(n)}}{h}\right) \right\} (=: \widehat{F}_{\infty}^{[R]}(X_{(n)})).
\end{align*}
The equation
$$\widehat{F}_{\infty}^{[R]}(X_{(n)}) < \frac{n}{n+1} < \widehat{F}_{X_{(n)}}^{[R]}(X_{(n)})$$
does not always hold, but we can prove the following asymptotic properties under some assumptions. We can get the following result in a similar manner.

\begin{cor} 
Let us assume that $f^{(1)}$ exists and is continuous and that $(7)$ holds. Then, if Assumption 1 holds and $h=o(n^{-1/4})$, we have
\begin{align*}
\widehat{\mbox{\boldmath $u$}}^{[R]} =& \mbox{\boldmath $u$}_0 + O_P(n^{-1/2} \mbox{\boldmath $1$}),\\
\widehat{F}_{\widehat{\mbox{\boldmath $u$}}}^{[R]} (x) =& F(x) + b_F^{[R]}(x) h^2 + O_P(n^{-1/2})\\ 
\intertext{and}
\widehat{f}_{\widehat{\mbox{\boldmath $u$}}}^{[R]} (x) =& f(x) + b_f^{[R]}(x) h + O_P((nh)^{-1/2}).
\end{align*}
\end{cor}

\section{Simulation study in Case 1}
\label{sect:Simulation study in Case 1}
We compared the proposed estimator numerically with the naive estimator and $X_{(n)}$-based estimator. Hall \& Park(2002) also discusses bias reduction of $X_{(n)}$ in case $f(u_0) > 0$; however, the modification does not affect the first-order asymptotics of the density estimation since $X_{(n)}$ is $n$ consistent in such case. Therefore, we did not take it into account here. Table 2 shows the averaged values of $ISE$, defined as
$$\int_{u_0 -h}^{\infty} \Bigm(\widetilde{f}(x) -f(x) \Bigm)^2 dx$$
of the density estimator $\widetilde{f}$ in the boundary region $\{ x > u_0-h \}$. In the table, $Naive$ denotes the $ISE$ of the standard kernel density estimator $\widehat{f}$. $\widehat{u}^{[BK]}$ denotes the $ISE$ of $\widehat{f}_{\widehat{u}}^{[BK]}$ and $X_{(n)}^{[BK]}$ denotes $X_{(n)}$-based estimator whose upper bound of support is given by $X_{(n)}$ using the boundary kernel method. $\widehat{u}^{[R]}$ denotes $\widehat{f}_{\widehat{u}}^{[R]}$ and $X_{(n)}^{[R]}$ is the $X_{(n)}$-based estimator using the boundary kernel method. The number of repetitions was $N=10000$ for all cases. The kernel functions were the Epanechnikov, and the bandwidths were the same and chosen by cross-validation as to which made `Naive' asymptotically best. $\widehat{f}_{\widehat{u}}^{[R]}$, $\widehat{f}$, and $\widehat{f}_{X_{(n)}}^{[R]}$ show the $ISE$ in the same way. The kernel function was a Gaussian.

\begin{table}
\caption{ISE of the kernel estimators}
$$\begin{tabu}[c]{c||c|c|c|c|c}
\hline

~~Beta(1,1)~~ & ~~Naive~~ & \widehat{u}^{[BK]} & X_{(n)}^{[BK]} & \widehat{u}^{[R]} & X_{(n)}^{[R]} \\ \hline

n=30 & .04414 & .06093 & .11531 & .04626 & .04646 \\

n=50 & .03823 & .03736 & .08202 & .03444 & .03452 \\

n=100  & .03161  & .02057 & .04270 & .02369 & .02374 \\

n=300  & .02340 & .00770 & .01732 & .01454  & .01457 \\ \hline \hline

Beta(3,1) & ~~Naive~~ & \widehat{u}^{[BK]} & X_{(n)}^{[BK]} & \widehat{u}^{[R]} & X_{(n)}^{[R]} \\ \hline

n=30 & .23968 & .18527 & .34458 & .21659 & .21811 \\

n=50 & .21025 & .11735 & .24758 & .17229 & .17299 \\

n=100  & .17250  & .05956 & .13326 & .12710  & .12752 \\

n=300  & .11425 & .02162 & .06087 & .07695  & .07726 \\ \hline \hline

Beta(3,3) & ~~Naive~~ & \widehat{u}^{[BK]} & X_{(n)}^{[BK]} & \widehat{u}^{[R]} & X_{(n)}^{[R]} \\ \hline

n=30 & .00596 & .02917 & .00834 & .00256 & .00256 \\

n=50 & .00359 & .01594 & .00308 & .00156 & .00156 \\

n=100  & .00173  & .00520 & .00152 & .00079  & .00079 \\

n=300  & .00058 & .00185 & .00041 & .00027  & .00027 \\

\hline
\end{tabu}$$
\end{table}

\begin{remark}
$f_{Beta(1,1)}(1) >0$, $f_{Beta(3,1)}(1) >0$, $f_{Beta(3,3)}(1) =0$ hold, where $f_{Beta(p,q)}$ denotes the density function of the beta distribution with parameters $(p,q)$. Note that if $f(1)=f'(1)=0$, $Bias[\widehat{f}(1)]=O(h^2)$ holds. Table 3 shows the boundary bias of the kernel estimators.
\end{remark}

\begin{table}
\caption{Boundary bias of the kernel estimators}
$$\begin{tabu}[c]{c|c|c|c|c|c}\hline 
Bias & ~~~\widehat{f}(1)~~~ & ~{\widehat{f}_{\widehat{u}}^{[BK]}}(1)~ & ~\widehat{f}_{X_{(n)}}^{[BK]}(1)~ & ~{\widehat{f}_{\widehat{u}}^{[R]}}(1)~ & ~\widehat{f}_{X_{(n)}}^{[R]}(1)~ \\ \hline
f(1) >0 & ~~~O(1)~~~ & \multicolumn{4}{c}{O(h)} \\ \hline
~f(1) =f'(1) =0~ & O(h^2) & \multicolumn{4}{c}{O(h^2)} \\ \hline
\end{tabu}$$
\end{table}

We can see from these numerical studies that when $f(1) >0$, $\widehat{u}$ converges more slowly than $X_{(n)}$, which has $n$ consistency; however, the proposed $\widehat{f}_{\widehat{u}}^{[\cdot]}$ is mostly better than $\widehat{f}_{X_{(n)}}^{[\cdot]}$. Although the boundary effect is not so great when $f(1)=0$, the proposed estimator is comparable with the other estimators. From the above results, we can claim that the proposed method is better in terms of both the theoretical and numerical local loss, at least when the support seems to be compact. Compared with the reflection method, the boundary kernel method seems to be numerically superior, especially when $f(1)>0$. 

\section{Extension to a simple multivariate case}
\label{Extension to a simple multivariate case} 
\subsection{Estimating multivariate joint density with unknown but simple form of support} 
We want to apply the new method to improve kernel joint probability density estimation, kernel regression and so on; however, it is impossible to apply it to multivariate cases directly. This is because the support of a multivariate joint density is given by a much more general formula and the problem completely changes. In addition, most of the boundary bias reduction methods have only been thoroughly investigated in simple cases, such as $[l_1,\infty) \times \cdots \times [l_d,\infty)$, $[0,1]^d$ and so on. 

Let $\mbox{\boldmath $X$}_1,\mbox{\boldmath $X$}_2,\cdots,\mbox{\boldmath $X$}_n$ be independently and identically distributed ($i.i.d.$) random variables with distribution function $G$ and density function $g$. $\mbox{\boldmath $X$}_i = (X_{1,i}, \cdots,X_{d,i})^T$ denotes $d$-dimensional value, and hereafter, we will assume that the support of $g$ is given by an unknown hyper-rectangle $\Xi$.
\begin{assumption}
There is bounded $l_{1,0}< u_{1,0}, \cdots ,l_{d,0} <u_{d,0}$ which satisfies
$$supp(g) = \Xi = \xi_1 \times \cdots \times \xi_d,$$
where $\xi_j$ is a bounded and open, half-open, or closed interval from $l_{j,0}$ to $u_{j,0}$.  
\end{assumption}

We can treat the marginal distributions in one dimension, even those whose support is unknown, so let us construct the joint probability density estimator, combining the marginals with the copula $C$ defined as follows:
\begin{eqnarray*}
C({\accentset{\diamond}{\mbox{\boldmath $G$}}}(\mbox{\boldmath $x$})) = G(\mbox{\boldmath $x$}),
\end{eqnarray*}
where $\mbox{\boldmath $x$} = (x_1,\cdots,x_d)^T$, ${\accentset{\diamond}{\mbox{\boldmath $G$}}}(\mbox{\boldmath $x$}) = (G_1(x_1),\cdots,G_d(x_d))^T$ and $G_j$ is the marginal distribution of the $j$-th component. The following equivalent form of the copula is known:
\begin{eqnarray*}
C(\mbox{\boldmath $v$})= G(G_1^{-1}(v_1),\cdots,G_d^{-1}(v_d)),
\end{eqnarray*}
where $\mbox{\boldmath $v$} = (v_1, \cdots, v_d)^T$ and the existence and uniqueness is known as Sklar's theorem. The copula function describes the dependence structure between the marginal variables. Scaillet and Fermanian (2002) proposed the following estimator in the bivariate case ($d=2$):
\begin{eqnarray}
\widetilde{C}(\mbox{\boldmath $v$}) = \frac{1}{n} \sum_{i=1}^n \prod_{j=1}^2 W\left( \frac{v_j -\widehat{G}_j(X_{j,i})}{h}\right),
\end{eqnarray}
where $\widehat{G}_j$ is the kernel marginal distribution estimator. In addition, Chen \& Huang (2007) discusses reduction of the boundary bias of the copula estimator; however, we can not apply it to the case of $d \ge 3$ because the boundary problem gets rather complicated. Instead, using an appropriate estimator $\widehat{C}_{n,h}$ and the idea of copula, we define the following general form of the joint probability density estimator free from the boundary bias
\begin{eqnarray}
&\widehat{G}_{\widehat{\Xi}}^{\dagger}(\mbox{\boldmath $x$}) := \widehat{C}_{n,h}\left( \widehat{\accentset{\diamond}{\mbox{\boldmath $G$}}}_{\widehat{\mbox{\boldmath $u$}}_1, \cdots, \widehat{\mbox{\boldmath $u$}}_d}^{\dagger}(\mbox{\boldmath $x$}) \right),
\end{eqnarray}
where
$$\widehat{\accentset{\diamond}{\mbox{\boldmath $G$}}}_{\widehat{\mbox{\boldmath $u$}}_1, \cdots, \widehat{\mbox{\boldmath $u$}}_d}^{\dagger}(\mbox{\boldmath $x$}) = \left( \widehat{G}_{1,\widehat{\mbox{\boldmath $u$}}_1}^{\dagger}(x_1), \cdots, \widehat{G}_{d,\widehat{\mbox{\boldmath $u$}}_d}^{\dagger}(x_d)\right)^T$$
$\widehat{G}_{j,\widehat{\mbox{\boldmath $u$}}_j}^{\dagger}$ is the proposed estimator of the marginal distribution $G_j$ as given in section 3 and $\widehat{C}_{n,h}$ satisfies the following assumption.

\begin{assumption}
For every $l_1=1,2$ and $l_2=1,2$, $\widehat{C}_{n,h}^{(l_1,l_2)}$ exists and is continuous. In addition, all of the following hold for any $\mbox{\boldmath $x$} \in \mathbb{R}^d$:
\begin{align*}
&E[\widehat{C}_{n,h}\left(\accentset{\diamond}{\mbox{\boldmath $G$}}(\mbox{\boldmath $x$})\right) ] = G(\mbox{\boldmath $x$}) + b_{G,\Xi}^{\dagger}(\mbox{\boldmath $x$},n,h), ~~~ V[\widehat{C}_{n,h}\left(\accentset{\diamond}{\mbox{\boldmath $G$}}(\mbox{\boldmath $x$})\right) ] = O(n^{-1}), \\
&E\left[ \left(\prod_{j=1}^d \frac{\partial}{\partial x_j} \right) \widehat{C}_{n,h}\left(\accentset{\diamond}{\mbox{\boldmath $G$}}(\mbox{\boldmath $x$})\right) \right] = g(\mbox{\boldmath $x$}) + b_{g,\Xi}^{\dagger}(\mbox{\boldmath $x$},n,h) \\ \intertext{and}
&V\left[\left(\prod_{j=1}^d \frac{\partial}{\partial x_j} \right) \widehat{C}_{n,h}\left(\accentset{\diamond}{\mbox{\boldmath $G$}}(\mbox{\boldmath $x$})\right) \right] = O((nh^d)^{-1}).
\end{align*}
\end{assumption} 

\begin{remark}
The following product-type of the boundary kernel and the reflection estimators of the marginals satisfy  Assumption 6 and give an appropriate distribution estimation,
$$\widehat{C}_{n,h}\left(\accentset{\diamond}{\mbox{\boldmath $G$}}(\mbox{\boldmath $x$})\right) = \frac{1}{n} \sum_{i=1}^n \prod_{j=1}^d W_{\mbox{\boldmath $u$}_{j,0}}^{[\cdot]} (x_j, X_{j,i},n,h),$$
where $n^{-1} \sum_{i=1}^n W_{\mbox{\boldmath $u$}_{j,0}}^{[\cdot]} (x_j, X_{j,i},n,h)= \widehat{G}_{\mbox{\boldmath $u$}_{j,0}}^{[\cdot]}(x_j)$ is the marginal distribution estimator of $G_j$ and `$\cdot$' is `$BK$' or `$R$'.
\end{remark}

The critical part of the modification $(14)$ is changing the unknown support (hypercube) to $[0,1]^d$. Under the assumption of support, we can remove the boundary bias.

\begin{theorem} Let us assume that for every marginal distribution $G_j$ $(j=1, \cdots, d)$, Assumptions 1 - 4 hold. In addition, we assume that for every $j$,
\begin{eqnarray}
\sqrt{n}(\mbox{\boldmath $X$}_{j,(1,n)} -\mbox{\boldmath $u$}_{j,0}) \xrightarrow{p} \mbox{\boldmath $0$}
\end{eqnarray}
holds, where $\mbox{\boldmath $X$}_{j,(1,n)} = (X_{j,(1)},X_{j,(n)})^T$ and $\mbox{\boldmath $u$}_{j,0} = (l_{j,0}, u_{j,0})^T$. Under Assumptions 5 and 6, we have
\begin{align}
\widehat{G}_{\widehat{\Xi}}^{\dagger} (\mbox{\boldmath $x$}) =& G(\mbox{\boldmath $x$}) + B_{G}^{\dagger}(\mbox{\boldmath $x$},n,h) + O_P(n^{-1/2}) \\ \intertext{and}
\widehat{g}_{\widehat{\Xi}}^{\dagger} (\mbox{\boldmath $x$}) =& g(\mbox{\boldmath $x$}) + B_{g}^{\dagger}(\mbox{\boldmath $x$},n,h) + O_P((nh^d)^{-1/2}),
\end{align}
where 
$$B_{G}^{\dagger}(\mbox{\boldmath $x$},n,h) = b_{G,\Xi}^{\dagger}(\mbox{\boldmath $x$},n,h) + \sum_{j=1}^{d} b_{G_j}^{\dagger}(\mbox{\boldmath $x$},n,h),$$
$b_{G_j}^{\dagger}(\mbox{\boldmath $x$},n,h)$ is the bias of the marginal distribution estimation $\widehat{G}_{{\mbox{\boldmath $u$}}_j,j}^{\dagger}(\mbox{\boldmath $x$})$ and $B_{g}^{\dagger}(\mbox{\boldmath $x$},n,h)$ is defined similarly.
\end{theorem}
{\it Proof.} See the appendices.

Table 4 summarizes the multivariate version of the proposed method.

\begin{table}
\caption{Flow of the proposed method in the multivariate case}
\begin{tabular}{l}
\hline \vspace{5pt}
0. Assume that the support of the density is of the form, \\ \vspace{5pt}
~~~~~ $supp(g)=\Xi=\xi_1 \times \cdots \times \xi_d$, \\ \vspace{5pt}
~~~~~ where $\xi_j$ is the interval from $l_{j,0}$ to $u_{j,0}$. \\ \vspace{3pt}
~~~ Check the assumptions of Theorem 3. \\ \vspace{3pt}
1. Select a boundary bias reduction method, \\ \vspace{3pt}
~~~~~ and choose the kernel and bandwidth as follows.\\ \vspace{3pt}
~~~ $\widehat{G}_{\mbox{\boldmath $u$}_j,j}^{\dagger}$ and $\widehat{g}_{\mbox{\boldmath $u$}_j,j}^{\dagger}$ denote the marginal distribution and density estimator, \\ \vspace{3pt}
~~~~~ free from the boundary bias caused by the bounded support $\mbox{\boldmath $u$}_j =[l_j,u_j]$.\\ \vspace{5pt}
2. Solve the equation for $\mbox{\boldmath $u$}_j$ for every $j$ : \\ \vspace{5pt}
\hspace{100pt} $\widehat{{\mbox{\boldmath $G$}}}_{\mbox{\boldmath $u$}_j,j}^{\dagger}(\mbox{\boldmath $X$}_{j,(1,n)}) - \mbox{\boldmath $c$}_n =\mbox{\boldmath $0$}$ \\
~~~~~ (see Section 3). \\ \vspace{3pt}
3. Set the solution as $\widehat{\mbox{\boldmath $u$}}_j$, \\ \vspace{3pt}
~~~~~ and output the boundary-adjusted estimator $\widehat{G}_{\widehat{\Xi}}^{\dagger}$ and $\widehat{g}_{\widehat{\Xi}}^{\dagger}$,\\ \vspace{3pt}
~~~~~ using a combining nonparametric copula estimator $\widehat{C}_{n,h}$ (see (14)).\\
\hline
\end{tabular}
\end{table}

\noindent\textbf{Acknowledgement}

The author thanks Prof. Y. Maesono., Faculty of Mathematics, Kyushu University for his valuable comments.

\section{Appendices: Some proofs}
\label{appendices}
\noindent
\textbf{Proof of Theorem 1}\medskip 
Under the assumptions of Theorem 1, we can prove that the following equation holds by using the asymptotic expansions:
$$\mbox{\boldmath $\Psi$}_{\mbox{\boldmath $u$},n}(\mbox{\boldmath $X$}_{(1,n)}) =\mbox{\boldmath $\Psi$}_{\mbox{\boldmath $u$},n}(\mbox{\boldmath $u$}_0) +o_P(n^{-1/2} \mbox{\boldmath $1$}).$$
Therefore, $\mbox{\boldmath $u$}=\widehat{\mbox{\boldmath $u$}}$ can be viewed as an $M$ estimator, and we can use asymptotic theory. Using assumptions $(i)$ and $(iv)$ of Assumption 4, we find that 
$$\|\mbox{\boldmath $\Psi$}_{\widehat{\mbox{\boldmath $u$}},n}(\mbox{\boldmath $u$}_0) -\mbox{\boldmath $\Psi$}_{\mbox{\boldmath $u$}_0,n}(\mbox{\boldmath $u$}_0)\| = o_P(1),$$
and
$$\|\mbox{\boldmath $\Psi$}_{\widehat{\mbox{\boldmath $u$}},n}(\mbox{\boldmath $u$}_0) -\mbox{\boldmath $F$}(\mbox{\boldmath $u$}_0)\| = o_P(1).$$
Combining them and using assumption $(ii)$ of Assumption 4 for any $\eta>0$, we have
\begin{align*}
P\left[ \|\widehat{\mbox{\boldmath $u$}} - \mbox{\boldmath $u$}_0 \| >  \eta \right] <& (1-\delta_n) P[\|\mbox{\boldmath $\Psi$}_{\mbox{\boldmath $u$}_0,n}(\mbox{\boldmath $u$}_0) - \mbox{\boldmath $\Psi$}_{\widehat{\mbox{\boldmath $u$}},n}(\mbox{\boldmath $u$}_0)\| > \kappa_{\eta,\delta,n}] + \delta_n \\
=&P[ \|\mbox{\boldmath $\Psi$}_{\mbox{\boldmath $u$}_0,n}(\mbox{\boldmath $u$}_0) - \mbox{\boldmath $\Psi$}_{\widehat{\mbox{\boldmath $u$}},n}(\mbox{\boldmath $u$}_0) \| > \kappa_{\eta,\delta,n}] + o(1) \\
\to& 0.
\end{align*}
Next, let us expand $\widehat{u}$ around $u_0$ in the estimating equation,
\begin{align*}
\mbox{\boldmath $0$} =& \widehat{\mbox{\boldmath $F$}}_{\widehat{\mbox{\boldmath $u$}}}^{\dagger}(\mbox{\boldmath $X$}_{(1,n)}) -\mbox{\boldmath $c$}_n \\
=& \mbox{\boldmath $\Psi$}_{\widehat{\mbox{\boldmath $u$}},n}(\mbox{\boldmath $u$}_0) -\mbox{\boldmath $c$}_n +o_P(n^{-1/2} \mbox{\boldmath $1$})\\
=& \mbox{\boldmath $\Psi$}_{\mbox{\boldmath $u$}_0,n}(\mbox{\boldmath $u$}_0) +  \left(\frac{\partial}{\partial \mbox{\boldmath $u$}} \mbox{\boldmath $\Psi$}_{\widetilde{\mbox{\boldmath $u$}},n}(\mbox{\boldmath $u$}_0)\right) (\widehat{\mbox{\boldmath $u$}} -\mbox{\boldmath $u$}_0) -\mbox{\boldmath $c$}_n +o_P(n^{-1/2} \mbox{\boldmath $1$}), 
\end{align*}
where $\widetilde{\mbox{\boldmath $u$}}$ is a random variable between $\widehat{\mbox{\boldmath $u$}}$ and $\mbox{\boldmath $u$}_0$. Accordingly, we find that
\begin{align*}
\widehat{\mbox{\boldmath $u$}} -\mbox{\boldmath $u$}_0 =& -\left[\frac{\partial}{\partial \mbox{\boldmath $u$}} \mbox{\boldmath $\Psi$}_{\widetilde{\mbox{\boldmath $u$}},n}(\mbox{\boldmath $u$}_0) \right]^{-1} (\mbox{\boldmath $\Psi$}_{\mbox{\boldmath $u$}_0,n}(\mbox{\boldmath $u$}_0) -\mbox{\boldmath $c$}_n +o_P(n^{-1/2} \mbox{\boldmath $1$}) ), \\
\end{align*}
where
$$\frac{\partial}{\partial \mbox{\boldmath $u$}} \mbox{\boldmath $\Psi$}_{\widetilde{\mbox{\boldmath $u$}},n}(\mbox{\boldmath $u$}_0) = E\left[\frac{\partial}{\partial \mbox{\boldmath $u$}} \mbox{\boldmath $\Psi$}_{\mbox{\boldmath $u$}_n^*,n}(\mbox{\boldmath $u$}_0) \right] + o_P(\mbox{\boldmath $1$} \otimes\mbox{\boldmath $1$})$$
holds. Therefore, we obtain
\begin{align*}
\widehat{\mbox{\boldmath $u$}} -\mbox{\boldmath $u$}_0 =& -E\left[\frac{\partial}{\partial \mbox{\boldmath $u$}} \mbox{\boldmath $\Psi$}_{\mbox{\boldmath $u$}_n^*,n}(\mbox{\boldmath $u$}_0) \right]^{-1} (\mbox{\boldmath $\Psi$}_{\mbox{\boldmath $u$}_0,n}(\mbox{\boldmath $u$}_0) -\mbox{\boldmath $c$}_n) +o_P(n^{-1/2} \mbox{\boldmath $1$}),
\end{align*}
and thus
$$\widehat{\mbox{\boldmath $u$}} -\mbox{\boldmath $u$}_0 = O_P(\mbox{\boldmath $\Psi$}_{\mbox{\boldmath $u$}_0,n}(\mbox{\boldmath $u$}_0) -\mbox{\boldmath $c$}_n) = O_P(n^{-1/2} \mbox{\boldmath $1$}).$$
Next, expanding $\widehat{\mbox{\boldmath $u$}}$ around $\mbox{\boldmath $u$}_0$, we have
\begin{align*}
\widehat{F}_{\widehat{\mbox{\boldmath $u$}}}^{\dagger}(x) =& \Psi_{\widehat{\mbox{\boldmath $u$}},n}(x) +R_{\widehat{\mbox{\boldmath $u$}},n}\\
=& \Psi_{\mbox{\boldmath $u$}_0,n}(x) + R_{\mbox{\boldmath $u$}_0,n} + \frac{1}{n} \sum_{i=1}^n \left(\frac{\partial}{\partial \mbox{\boldmath $u$}}\psi_{i,\mbox{\boldmath $u$}_0\,n} \left(x\right) \right)^T  (\widehat{\mbox{\boldmath $u$}}- \mbox{\boldmath $u$}_0) \\
& + O_P(\|\widehat{\mbox{\boldmath $u$}}- \mbox{\boldmath $u$}_0\|^2) + o_P(n^{-1}).
\end{align*}
The expectation is given by
\begin{align*}
E[\widehat{F}_{\widehat{\mbox{\boldmath $u$}}}^{\dagger}(x)] =& F(x) + b_F^{\dagger}(x,n,h) + O(n^{-1/2}).
\end{align*}

To derive the variance term, we perform another asymptotic expansion as follows:
\begin{align*}
[\widehat{F}_{\widehat{\mbox{\boldmath $u$}}}^{\dagger}(x)]^2 =& [\Psi_{\widehat{\mbox{\boldmath $u$}},n}(x) +R_{\widehat{\mbox{\boldmath $u$}},n}]^2 \\
=& \left[\Psi_{\mbox{\boldmath $u$}_0,n}(x) + R_{\mbox{\boldmath $u$}_0,n} + \frac{1}{n} \sum_{i=1}^n \left(\frac{\partial}{\partial \mbox{\boldmath $u$}}\psi_{i,\mbox{\boldmath $u$}_0,n} \left(x\right) \right)^T (\widehat{\mbox{\boldmath $u$}}- \mbox{\boldmath $u$}_0) \right]^2 + O_P(n^{-1})\\
=& \Psi_{\mbox{\boldmath $u$}_0,n}(x)\left[\Psi_{\mbox{\boldmath $u$}_0,n}(x) + R_{\mbox{\boldmath $u$}_0,n} + \frac{1}{n} \sum_{i=1}^n \left(\frac{\partial}{\partial \mbox{\boldmath $u$}}\psi_{\mbox{\boldmath $u$}_0,n} \left(x\right) \right)^T (\widehat{\mbox{\boldmath $u$}}- \mbox{\boldmath $u$}_0) \right] \\
&+ O_P(n^{-1}).
\end{align*}
Then, we have the following four results:
\begin{align*}
&E[(\widehat{F}_{\widehat{\mbox{\boldmath $u$}}}^{\dagger}(x))^2] \\
=& E\left[\Psi_{\mbox{\boldmath $u$}_0,n}(x)\left[\Psi_{\mbox{\boldmath $u$}_0,n}(x) + R_{\mbox{\boldmath $u$}_0,n} + \frac{1}{n} \sum_{i=1}^n \left(\frac{\partial}{\partial \mbox{\boldmath $u$}}\psi_{i,\mbox{\boldmath $u$}_0,n} \left(x\right) \right)^T (\widehat{\mbox{\boldmath $u$}}- \mbox{\boldmath $u$}_0) \right] \right] \\
&~~~ + O(n^{-1}), \\
&E\left[(\Psi_{\mbox{\boldmath $u$}_0,n}(x))^2 \right] -E\left[\Psi_{\mbox{\boldmath $u$}_0,n}(x) \right]^2 = O(n^{-1}),\\
&E\left[\Psi_{\mbox{\boldmath $u$}_0,n}(x) R_{\mbox{\boldmath $u$}_0,n} \right] - E\left[\Psi_{\mbox{\boldmath $u$}_0,n}(x)\right] E\left[R_{\mbox{\boldmath $u$}_0,n} \right] = O(n^{-1}),\\ \intertext{and}
&E\left[\Psi_{\mbox{\boldmath $u$}_0,n}(x) \frac{1}{n} \sum_{i=1}^n \left(\frac{\partial}{\partial \mbox{\boldmath $u$}}\psi_{i,\mbox{\boldmath $u$}_0,n} \left(x\right) \right)^T (\widehat{\mbox{\boldmath $u$}}- \mbox{\boldmath $u$}_0) \right]\\
& -E\left[\Psi_{\mbox{\boldmath $u$}_0,n}(x)\right] E\left[\frac{1}{n} \sum_{i=1}^n \left(\frac{\partial}{\partial \mbox{\boldmath $u$}}\psi_{i,\mbox{\boldmath $u$}_0,n} \left(x\right) \right)^T (\widehat{\mbox{\boldmath $u$}}- \mbox{\boldmath $u$}_0) \right] = O(n^{-1}).
\end{align*}
From the above results, we can see that
$$V[\widehat{F}_{\widehat{\mbox{\boldmath $u$}}}^{\dagger}(x)] = O(n^{-1}),$$
and that
$$\widehat{F}_{\widehat{\mbox{\boldmath $u$}}}^{\dagger} (x)= F(x) + b_F^{\dagger}(x,n,h) + O_P(n^{-1/2}).$$
In the same way, we can prove that
$$\widehat{f}_{\widehat{\mbox{\boldmath $u$}}}^{\dagger} (x)= f(x) + b_f^{\dagger}(x,n,h) + O_P((nh)^{-1/2}).$$

\noindent
\textbf{Proof of Theorem 2}
From Theorem 1, the following asymptotic expansion holds under the assumptions:
\begin{align*}
&\widehat{F}_{\widehat{\Xi}}^{\dagger}(\mbox{\boldmath $x$})\\
=& \widehat{C}_{n,h}\left( \accentset{\diamond}{\mbox{\boldmath $G$}}(\mbox{\boldmath $x$}) \right) + O\left(\sum_{j=1}^{d} b_{F_j}^{\dagger}(x_j,n,h)\right) + O_P(n^{-1/2}).
\end{align*}
Then, it is easy to see that Theorem 2 follows.

\end{document}